\documentclass{article}

\catcode`\@=11 \@addtoreset{equation}{section}

\catcode`\@=12

\newtheorem{Theorem}{Theorem}[section]
\newtheorem{Proposition}{Proposition}[section]
\newtheorem{Lemma}{Lemma}[section]
\newtheorem{Corollary}{Corollary}[section]

\newcommand{\bTheorem}[1]{
\begin{Theorem} \label{T#1} }
\newcommand{\eT}{\end{Theorem}}

\newcommand{\bProposition}[1]{
\begin{Proposition} \label{P#1}}
\newcommand{\eP}{\end{Proposition}}

\newcommand{\bLemma}[1]{
\begin{Lemma} \label{L#1} }
\newcommand{\eL}{\end{Lemma}}

\newcommand{\bCorollary}[1]{
\begin{Corollary} \label{C#1} }
\newcommand{\eC}{\end{Corollary}}

\newcommand{\bFormula}[1]{
\begin{equation} \label{#1}}
\newcommand{\eF}{\end{equation}}

\newcommand{\Ov}[1]{\overline{#1}}
\newcommand{\DC}{C^\infty_c}

\newcommand{\vr}{\varrho}

\newcommand{\vt}{\vartheta}
\newcommand{\vu}{\vc{u}}
\newcommand{\vc}[1]{{\bf #1}}

\newcommand{\Div}{{\rm div}_x}
\newcommand{\Grad}{\nabla_x}

\newcommand{\tn}[1]{\mbox {\F #1}}
\newcommand{\dx}{{\rm d} {x}}
\newcommand{\dt}{{\rm d} t }

\newcommand{\intO}[1]{\int_{\Omega} #1 \ \dx}

\font\F=msbm10 scaled 1000

\date{}
\makeindex
\begin{document}



\title{Relative entropies, suitable weak solutions, and weak-strong uniqueness
for the compressible Navier-Stokes system}
\author{Eduard Feireisl\thanks{The work of E.F. was supported by Grant 201/09/0917
of GA \v CR as a part of the general research programme of the
Academy of Sciences of the Czech Republic, Institutional Research
Plan AV0Z10190503.} \and Bum Ja Jin\thanks{The work of B.J.J. was supported by Basic
Science Research Program through the National Research Foundation
of Korea funded by the Ministry of Education, Science and Technology(2011-0007701).}\and Anton\' \i n
Novotn\' y
\thanks{The work was completed during the stay of A.N. at the Mathematical Institute
in Prague supported by the general research programme of the
Academy of Sciences of the Czech Republic, Institutional Research
Plan AV0Z10190503. A part of this work was effectuated during the
stay of A.N. at the National Mokpo University  under the financial
support of Science Research Program through the National Research
Foundation of Korea funded by the Ministry of Education, Science
and Technology(2011-0007701).}} \maketitle

\bigskip

\centerline{Institute of Mathematics of the Academy of Sciences of the Czech Republic}
\centerline{\v Zitn\' a 25, 115 67 Praha 1, Czech Republic}

\bigskip

\centerline{Department of Mathematics Education, Mokpo National University}
\centerline{Muan 534-729, South Korea}

\bigskip

\centerline{IMATH Universit\' e du Sud Toulon-Var} \centerline{BP
132, 83957 La Garde, France}

\bigskip

\bigskip

\begin{abstract}

We introduce the notion of \emph{relative entropy} for the weak solutions to the compressible Navier-Stokes system. In particular, we show that any
finite energy weak solution satisfies a relative entropy inequality with respect to any couple of smooth functions satisfying relevant boundary
conditions. As a corollary, we establish the weak-strong uniqueness property in the class of finite energy weak solutions, extending thus the
classical result of Prodi and Serrin to the class of compressible fluid flows.
\end{abstract}

\section{Introduction}
\label{i}

The method of \emph{relative entropy} has been successfully applied to partial differential equations of different types. \emph{Relative entropies}
are non-negative quantities that provide a kind of distance between two solutions of the same problem, one of which typically enjoys some extra
regularity properties. Carillo et al. \cite{CaJuMaToUn} exploited entropy dissipation, expressed by means of the relative entropy with respect to
a stationary solution, in order to analyze the long-time behavior of certain
quasilinear parabolic equations. Saint-Raymond \cite{SaRay} uses the relative entropy method to study the incompressible Euler limit of the Boltzmann
equation. Other applications of the method can be found in Grenier \cite{Greni}, Masmoudi \cite{MAS5}, Ukai \cite{Uka}, Wang and Jiang \cite{WanJia},
among others.

Germain \cite{Ger} introduced a class of (weak) solutions to the compressible Navier-Stokes system satisfying a relative entropy inequality with
respect to a (hypothetical) strong solution of the same problem, and established the weak-strong uniqueness property within this class. Unfortunately,
\emph{existence} of solutions belonging to this class, where, in particular, the density possesses a spatial gradient in a suitable Lebesgue space, is
not known. In \cite{FENOSU}, we introduced
the concept of \emph{suitable weak solution} for the compressible Navier-Stokes system, satisfying a general relative entropy inequality with respect
to any sufficiently regular pair of functions. To be more specific,
consider the fluid density $\vr = \vr(t,x)$, together with the velocity field
$\vu = \vu(t,x)$, $t \in R$, $x \in \Omega \subset R^3$, the time evolution of which is governed by the
\emph{Navier-Stokes system}:
\bFormula{i1}
\partial_t \vr + \Div (\vr \vu) = 0,
\eF
\bFormula{i2}
\partial_t (\vr \vu) + \Div (\vr \vu \otimes \vu) + \Grad p(\vr) =
\Div \tn{S}(\Grad \vu) + \vr \vc{f}, \eF \bFormula{i3} \tn{S} =
\mu \Big( \Grad \vu + \Grad^t \vu - \frac{2}{3} \Div \vu \tn{I}
\Big) + \eta \Div \vu \tn{I} , \ \mu > 0, \ \eta \geq 0, \eF
supplemented with suitable boundary conditions, say, \bFormula{i4}
\vu|_{\partial \Omega} = 0. \eF
If the domain $\Omega$ is unbounded, we prescribe the
far-field behavior: \bFormula{i4a} \vr \to \Ov{\vr}, \ \vu \to 0 \
\mbox{as}\ |x| \to \infty, \eF where $\overline\vr\ge 0$.

\emph{Relative entropy} ${\cal E}\Big( [\vr, \vu] \Big| [r ,
\vc{U}] \Big)$ with respect to $[r, \vc{U}]$ is defined as
\bFormula{i5} {\cal E} \Big( [\vr, \vu] \Big| [r, \vc{U}] \Big) =
\intO{ \left( \frac{1}{2} \vr |\vu - \vc{U}|^2 + H(\vr) - H'(r)
(\vr - r) - H(r) \right)}, \eF where \bFormula{i5a} H(\vr) = \vr
\int_{\Ov{\vr}}^\vr \frac{p(z)}{z^2} \ {\rm d}z.
\eF

Following \cite{FENOSU}, we say that $\vr$, $\vc{u}$ is a
\emph{suitable weak solution} to problem (\ref{i1} - \ref{i4a}) if
equations (\ref{i1}--\ref{i3}) are satisfied in a weak sense, and,
in addition to (\ref{i1} - \ref{i4a}), the following (relative)
energy inequality \bFormula{i6} {\cal E} \Big( [\vr, \vu] \Big|
[r, \vc{U}] \Big) (\tau) + \int_0^\tau \intO{ \Big( \tn{S} (\Grad
\vu) - \tn{S} (\Grad \vc{U}) \Big): \Big( \Grad \vu - \Grad \vc{U}
\Big) } \ \dt \eF
\[
\leq {\cal E} \Big([\vr_0, \vu_0] \Big| [r(0, \cdot), \vc{U}(0, \cdot)] \Big) +
\int_0^\tau {\cal R}(\vr, \vu, r, \vc{U} ) \ \dt
\]
holds for a.a. $\tau > 0$, where
\[
\vr_0 = \vr(0, \cdot), \ \vu_0 = \vu(0, \cdot),
\]
and the remainder ${\cal R}$ reads
\bFormula{i7} {\cal R}\left( \vr, \vu, r, \vc{U} \right) \equiv
\intO{  \vr \Big( \partial_t \vc{U} + \vu \Grad \vc{U} \Big) \cdot
(\vc{U} - \vu )}
\eF
\[
+ \intO{\tn{S}(\Grad \vc{U}):\Grad (\vc U- \vc{u})  }
+\intO{\vr\vc f \cdot(\vc u-\vc U)}
\]
\[
+ \intO{ \left( (r - \vr) \partial_t H'(r) + \Grad H'(r) \cdot
\left( r \vc{U} - \vr \vu \right) \right) }
- \intO{
\Div \vc{U} \Big( p(\vr) - p(r) \Big) }.
\]
Here, the functions $r$, $\vc{U}$ are arbitrary smooth, $r$
strictly positive, and $\vc{U}$ satisfying the no-slip boundary
conditions (\ref{i4}). It is easy to check that (\ref{i6}) is
satisfied as an equality as soon as the solution $\vr$, $\vu$ is
smooth enough.

As shown in \cite[Theorem 3.1]{FENOSU}, the Navier-Stokes system
(\ref{i1} - \ref{i4a}) admits global-in-time suitable weak
solutions for any finite energy initial data. Moreover, the
relative energy inequality (\ref{i6}) can be used to show that
suitable weak solutions comply with the weak-strong uniqueness
principle, meaning, a weak and strong solution emanating from the
same initial data coincide as long as the latter exists. This can
be seen by taking the strong solution as the ``test'' functions
$r$, $\vc{U}$ in the relative entropy inequality (\ref{i6}).
Besides, a number of other interesting properties of the suitable
weak solutions can be deduced, see \cite[Section 4]{FENOSU}.

For the particular choice $r = \Ov{\vr}$, $\vc{U} = 0$, the relative energy inequality (\ref{i6}) reduces to the standard \emph{energy inequality}
\bFormula{i8}
{\cal E} [\vr, \vu]  (\tau)
+ \int_0^\tau \intO{  \tn{S} (\Grad \vu):
\Grad \vu  } \ \dt  \leq {\cal E}  [\vr_0, \vu_0]  +
\int_0^\tau \intO{ \vr \vc{f} \cdot \vu } \ \dt \ \mbox{for a.a.}\
\tau > 0,
\eF
\[
{\cal E}[\vr, \vu] = \intO{ \left( \frac{1}{2} \vr |\vu|^2 +
H(\vr) - H'(\Ov{\vr}) \Big( \vr - \Ov{\vr}\Big) -H(\overline\vr)
\right) }.
\]
The weak solutions of the Navier-Stokes system satisfying, in addition, the energy inequality
(\ref{i8}) are usually termed \emph{finite energy weak solutions}, or, rather incorrectly,
turbulent solutions in the sense of Leray's original work \cite{LER}.

Our goal in this paper is to show that any finite energy weak
solution is in fact a suitable weak solution, in other words, the
standard energy inequality (\ref{i8}) implies the relative energy
inequality (\ref{i6}). In particular, the weak-strong uniqueness
property as well as other results shown in \cite{FENOSU} hold for
the seemingly larger class of finite energy solutions. This
observation extends easily to other types of boundary conditions
and to a large class of domains. This kind of result can be viewed
as an extension of the seminal work of Prodi \cite{PR} and Serrin
\cite{SERRIN1} (see also Germain \cite{Ger2} for more recent
results) to the compressible Navier-Stokes system. We provide an
ultimate answer to the weak-strong uniqueness problem intimately
related to the fundamental questions of the well-posedness for the
compressible Navier-Stokes equations addressed by several authors,
Desjardin \cite{DES2}, Germain \cite{Ger}, Hoff \cite{Hoff10},
\cite{Hoff9}, among others.

The paper is organized as follows. In Section \ref{w}, we provide an exact definition of finite
energy weak solutions and state the main result. Section
\ref{p} is devoted to the proof of the main theorem and
to possible extensions. Applications are discussed in Section \ref{a}.

\section{Main results}
\label{w} For the sake of simplicity, we assume that the pressure
$p = p(\vr)$ is a continuously differentiable function of the
density such that \bFormula{w1} p \in C[0,\infty) \cap
C^2(0,\infty), \ p(0) = 0,\ p'(\vr) > 0 \ \mbox{for all}\ \vr > 0,
\ \lim_{\vr \to \infty} \frac{ p'(\vr) } {\vr^{\gamma - 1}} = a >
0 \ \mbox{for a certain}\ \gamma > 3/2. \eF
Moreover, if $\overline\vr=0$, we suppose that $p$ becomes asymptotically small for
$\vr \to 0$ so that the function $H$ defined
in (\ref{i5a}) is finite for any $\vr>0$.

\subsection{Finite energy weak solutions to the Navier-Stokes system}
\label{fws}
{\bf Definition \ref{w}.1} \\
{\it We shall say that $\vr$, $\vu$ is a \emph{finite energy weak
solution} to the Navier-Stokes system (\ref{i1} - \ref{i4a})
emanating from the initial data $\vr_0$, $\vu_0$ if
\begin{itemize}
\item
\begin{equation}\label{d1}
\vr - \Ov{\vr} \in L^\infty(0,T;L^2 + L^\gamma(\Omega)),\; \vr
\geq 0\;\mbox{ a.a. in $(0,T) \times \Omega)$};
\end{equation}
\begin{equation}\label{d2}
\vu \in L^2(0,T;D^{1,2}_0(\Omega;R^3));
\end{equation}
\begin{equation}\label{d3}
\vr \vu \in L^\infty(0,T;L^2 + L^{2\gamma/(\gamma +
1)}(\Omega;R^3));
\end{equation}
\begin{equation}\label{d4}
p \in L^1_{\rm loc}([0,T] \times \Omega);
\end{equation}
\item $(\vr - \Ov{\vr}) \in C_{\rm weak}([0,T]; L^2 + L^\gamma(\Omega))$ and the integral identity
\bFormula{w2}
\intO{ \vr (\tau, \cdot) \varphi (\tau, \cdot) } -
\intO{ \vr_0 \varphi(0, \cdot) } =
\int_0^T \intO{ \Big( \vr \partial_t \varphi + \vr \vu \cdot \Grad \varphi \Big) } \ \dt
\eF
holds for any $\varphi \in \DC([0,T] \times \Ov{\Omega})$;
\item $\vr \vu \in C_{\rm weak}([0,T]; L^2 + L^{2\gamma/(\gamma + 1)}(\Omega;R^3))$
and the integral identity
\bFormula{w3}
\intO{ \vr \vu (\tau, \cdot) \cdot \varphi (\tau, \cdot) } -
\intO{ \vr_0 \vu_0 \cdot \varphi(0, \cdot) }
\eF
\[
= \int_0^T \intO{ \Big( \vr \vu \cdot \partial_t \varphi + (\vr \vu \otimes \vu): \Grad \varphi + p(\vr)
\Div \varphi - \tn{S} (\Grad \vu) : \Grad \varphi + \vr \vc{f} \cdot \varphi \Big) } \ \dt
\]
is satisfied for any
$\varphi \in \DC([0,T] \times \Omega ; R^3)$;
\item
the energy inequality \bFormula{w4} \intO{ \Big( \frac{1}{2} \vr
|\vu|^2 + H(\vr) - H'(\Ov{\vr}) (\vr - \Ov{\vr})-H(\overline\vr)
\Big) (\tau, \cdot) } + \int_0^\tau \intO{ \tn{S}(\Grad \vu) :
\Grad \vu } \ \dt \eF
\[
\leq \intO{ \Big( \frac{1}{2} \vr_0 |\vu_0|^2 + H(\vr_0) -
H'(\Ov{\vr}) (\vr_0 - \Ov{\vr}) -H(\overline\vr) \Big) } +
\int_0^T \intO{ \vr \vc{f} \cdot \vu } \ \dt
\]
holds for a.a. $\tau \in [0,T]$.
\end{itemize}
}

{\bf Remark \ref{w}.1} {\it We recall that the space
$D^{1,2}_0(\Omega)$ is defined as a completion of
$\DC(\Omega)$ with respect to the $L^2-$norm of the
gradient. In accordance with Sobolev's inequality,
\[
D^{1,2}_0(\Omega) \subset L^6(\Omega),
\]
see Galdi \cite{GAL}. }

\medskip

{\bf Remark \ref{w}.2} {\it In (\ref{w4}), we tacitly assume that
the initial data are chosen in such a way that the first integral
on the right hand side is finite.
}

\medskip

\subsection{Finite energy weak solutions satisfy the relative energy inequality}

Our main result reads as follows:

\bTheorem{w1} Let $\Omega \subset R^3$ be a domain.
Suppose that the pressure $p$ satisfies hypothesis (\ref{w1}),
\[
\vc{f} \in L^\infty(0,T; L^1 \cap L^\infty (\Omega; R^3)),
\]
and that  $\Ov{\vr} \ge 0$. Let $\vr$, $\vu$ be a finite energy
weak solution to the Navier-Stokes system (\ref{i1} - \ref{i4a})
in the sense specified in Section \ref{fws}.

Then $\vr$, $\vu$ satisfy the relative energy inequality (\ref{i6}) for any
$\vc{U} \in \DC([0,T] \times \Omega;R^3)$, and $r > 0$,
$r - \Ov{\vr} \in \DC([0,T] \times \Ov{\Omega})$.
\eT

The proof and several extensions of Theorem \ref{Tw1} are presented in Section \ref{p}. Applications will
be discussed in Section \ref{a}.

\section{Proof of the main result}
\label{p}

\subsection{Proof of Theorem \ref{Tw1}}

Take $\vc{U}$ as a test function in the momentum equation (\ref{w3}) to obtain
\bFormula{p1}
\intO{ \vr \vu (\tau, \cdot) \cdot \vc{U} (\tau, \cdot) } =
\intO{ \vr_0 \vu_0 \cdot \vc{U} (0, \cdot) }
\eF
\[
+ \int_0^\tau \intO{ \Big( \vr \vu \cdot \partial_t \vc{U} + (\vr \vu \otimes \vu): \Grad \vc{U} + p(\vr)
\Div \vc{U} - \tn{S} (\Grad \vu) : \Grad \vc{U} + \vr \vc{f} \cdot \vc{U} \Big) } \ \dt
\]
Similarly, we can use the scalar quantity $\frac{1}{2} |\vc{U}|^2$ as a test function in (\ref{w2}):
\bFormula{p2}
\intO{ \frac{1}{2} \vr (\tau, \cdot) |\vc{U}|^2 (\tau, \cdot) } =
\intO{ \frac{1}{2} \vr_0 |\vc{U} (0, \cdot)|^2 } +
\int_0^\tau \intO{ \Big( \vr \vc{U} \cdot \partial_t \vc{U} + \vr \vu \cdot \Grad \vc{U} \cdot \vc{U} \Big) } \ \dt .
\eF
Finally, we test (\ref{w2}) on $H'(r) - H'(\Ov{\vr})$ to get
\bFormula{p3}
\intO{  \vr (\tau, \cdot) \Big( H'(r) (\tau, \cdot) - H'(\Ov{\vr}) \Big) } =
\intO{  \vr_0 \Big( H'(r) (0, \cdot) - H'(\Ov{\vr}) \Big) }
\eF
\[
 +
\int_0^\tau \intO{ \Big( \vr \partial_t H'(r) + \vr \vu \cdot \Grad H'(r) \Big) } \ \dt .
\]

Summing up relations (\ref{p1} - \ref{p3}) with the energy inequality
(\ref{w4}), we infer that
\bFormula{p4}
\intO{ \left( \frac{1}{2} \vr |\vu - \vc{U}|^2 + H(\vr) -
\Big( H'(r) \vr - H'(\Ov{\vr}) \Ov{\vr} \Big) \right) (\tau, \cdot) }
\eF
\[
+ \int_0^\tau \intO{ \Big( \tn{S} (\Grad \vu) - \tn{S} (\Grad \vc{U}) \Big)
: \Big( \Grad \vu - \Grad \vc{U} \Big) } \ \dt
\]
\[
=\intO{ \left( \frac{1}{2} \vr_0 |\vu_0 - \vc{U} (0, \cdot) |^2 + H(\vr_0) -
\Big( H'(r(0, \cdot)) \vr_0 - H'(\Ov{\vr}) \Ov{\vr} \Big) \right)  }
\]
\[
+ \int_0^\tau \intO{ \vr \Big(\partial_t \vc{U} + \vr \vu \cdot \Grad \vc{U} \Big)\cdot (\vc{U} - \vu) } \ \dt
\]
\[
+ \int_0^\tau \intO{\tn{S}(\Grad \vc{U}):\Grad (\vc U- \vc{u})  }
+\intO{\vr\vc f \cdot(\vc u-\vc U)} \ \dt
\]
\[
-\int_0^\tau \intO{ \Big( \vr \partial_t H'(r) + \vr \vu \cdot \Grad H'(r) \Big) } \ \dt -  \int_0^\tau \intO{ p(\vr) \Div \vc{U} } \ \dt.
\]

Realizing that
\[
H'(r)r - H(r) - H'(\Ov{\vr}) \Ov{\vr} = p(r) - p(\Ov{\vr}),
\]
we compute
\[
\intO{ \Big( p(r) - p(\Ov{\vr}) \Big) (\tau, \cdot) } -
\intO{ \Big( p(r) - p(\Ov{\vr}) \Big) (0, \cdot) }
= \int_0^\tau \intO{ \partial_t p(r) } \ \dt ;
\]
whence, by virtue of the identity
\bFormula{p5}
\intO{ \left( r \partial_t H'(r) + r \Grad H'(r) \cdot \vc{U} + p(r)
\Div \vc{U} \right) } = \intO{ \partial_t p(r) },
\eF
relation (\ref{p4}) implies (\ref{i6}). Theorem \ref{Tw1} has been proved.
Note that (\ref{p5}) relies on the fact that $\vc{U} \cdot \vc{n}|_{\partial \Omega} = 0$.

\subsection{Possible extensions}

The conclusion of Theorem \ref{Tw1} can be extended in several
directions. Here, we shortly discuss the problem of an alternative choice of boundary conditions as well
as the class of admissible test functions $r$, $\vc{U}$.


\subsubsection{General slip boundary conditions with friction}

Similar result can be obtained provided the no-slip boundary
condition (\ref{i4}) is replaced by the  slip boundary
conditions with friction (Navier's boundary condition)  \bFormula{p6} \vu \cdot\vc n= 0 ,\;(\tn
S(\Grad \vc u)\vc n)_{\rm tan} +\beta \vc u_{\rm tan}=0\
\mbox{on} \ (0,T)\times
\partial \Omega, \eF
where $\beta\ge 0$ and $\vc v_{\rm tan}|_{\partial\Omega}=(\vc
v-(\vc v\cdot\vc n)\vc n)|_{\partial\Omega}$ denotes the tangential
componenet of a vector field $\vc{v}$ at the boundary. Note that the
so-called complete slip boundary conditions correspond to the particular
sitution $\beta=0$.

The definition of finite energy weak solutions is similar
to Section \ref{fws} with the following modifications:
\begin{itemize}
\item the spatial domain $\Omega$ possesses a Lipschitz boundary, where (\ref{d2}) is replaced by the requirement
$\vu \in L^2(0,T; D^{1,2}_n(\Omega;R^3))$, with
\[
D^{1,2}_n (\Omega;R^3) = \left\{ \vc{v} \in L^6_{\rm loc}(\Ov{\Omega};R^3) \ \Big| \
\Grad \vc{v} \in L^2(\Omega; R^{3 \times 3}), \ \vc{v} \cdot \vc{n}|_{\partial \Omega} = 0 \right\};
\]

\item the pressure satisfies
\begin{equation}\label{d4+}
p(\vr) \in L^1_{\rm loc}([0,T] \times \Ov{\Omega})
\end{equation}
instead of (\ref{d4});
\item the weak formulation of the momentum equation (\ref{w2}) has
to be replaced by
\bFormula{w2+} \int_0^\tau \intO{ \Big( \vr \vu \cdot \partial_t
\varphi + \vr (\vu \otimes \vu) : \Grad \varphi + p(\vr) \Div
\varphi \Big) } \ \dt \eF
\[
 -  \int_0^\tau \intO{ \tn{S}(\Grad \vu) : \Grad \varphi } \
\dt -\beta\int_0^\tau\int_{\partial\Omega}\vc u\cdot\varphi{\rm d S}\ \dt
\]
\[
 =-\int_0^\tau\intO{\vr\vc f\cdot\varphi}\ \dt+ \intO{
(\vr\vc{u})(\tau) \cdot \varphi (\tau, \cdot) } - \intO{ \vr_0
\vc{u}_0 \cdot \varphi (0, \cdot) }
\]
for all $\tau\in [0,T]$, for any test function $\varphi \in
\DC([0,T] \times \overline\Omega; R^3)$, $\vc \varphi\cdot\vc n=0$
on $[0,T]\times\partial\Omega$;
\item energy inequality (\ref{w3}) is replaced by
\bFormula{w3+} \intO{ \left( \frac{1}{2} \vr |\vu|^2 + H(\vr) -
H'(\Ov{\vr}) (\vr - \Ov{\vr})-H(\overline\vr) \Big) (\tau, \cdot)
\right)(\tau, \cdot) } \eF \[
 + \int_0^\tau \intO{ \tn{S}(\Grad \vu) :
\Grad \vu } \ \dt +\beta\int_0^\tau\int_{\partial\Omega}|\vc
u |^2{\rm d S}\ \dt \]
\[
\leq\int_0^\tau\intO{\vr\vc f\cdot\vc u}\ \dt+ \intO{ \left(
\frac{1}{2} \vr_0 |\vc{u}_0|^2 + H(\vr_0) - H'(\Ov{\vr}) (\vr_0 -
\Ov{\vr})-H(\overline\vr) \Big) (\tau, \cdot)\right) } \ \mbox{for
a.a.}\ \tau \in (0,T).
\]
\end{itemize}
%
In this case, the conclusion of Theorem \ref{Tw1} remains valid
for any couple $(r,\vc U)$ such that \bFormula{p6+}r-\overline\vr
\in C^\infty_c([0,T]\times\overline\Omega),\quad\vc{U} \in \DC
([0,T] \times \Ov{\Omega};R^3), \ \vc{U} \cdot \vc{n}|_{\partial
\Omega} = 0 \eF with the relative entropy inequality that reads
\bFormula{w5+}   {\cal E}([\vr,\vu]\Big| [r,\vc U])  (\tau, \cdot)
\eF
\[
 + \int_0^\tau \intO{ \left[
\tn{S} (\Grad \vu -\Grad \vc{U}) \right] : \Grad (\vc{u} - \vc{U})
} \ \dt + \beta\int_0^\tau\int_{\partial\Omega}|\vc u -\vc
U |^2{\rm d} S{\rm d}t
\]
\[
\leq {\cal E}([\vr_0,\vu_0]\Big| [r(0),\vc U(0)) (\tau ) +
\int_0^\tau {\cal R} \left( \vr, \vu, r, \vc{U} \right) \ \dt \
\mbox{for a.a.} \ \tau \in (0,T),
\]
where \bFormula{w6+} {\cal R}\left( \vr, \vu, r, \vc{U} \right)
=\intO{\vr\vc f\cdot(\vc u-\vc U)}
-\beta\int_0^\tau\int_{\partial\Omega}\vc U \cdot(\vc
u  -\vc U ){\rm d}S{\rm d}t \eF
\[
+ \intO{ \left( \vr \Big(
\partial_t \vc{U} + \vu\cdot \Grad \vc{U} \Big) \cdot (\vc{U} -
\vu ) -\tn{S}(\Grad \vc{U}) :\Grad (\vu - \vc{U}) \right) } \]
\[
+ \intO{ \left( (r - \vr) \partial_t H'(r) + \Grad H'(r) \cdot
\left( r \vc{U} - \vr \vu \right) - \Div \vc{U} \Big( p(\vr) -
p(r) \Big) \right) }.
\]

\subsubsection{Extending the admissible class of test functions}

Using density arguments we can extend considerably the class of
test functions $r$, $\vc{U}$ appearing in the relative energy
inequality (\ref{i6}) resp. (\ref{w5+}). Indeed:
\begin{itemize}
\item
For the left hand side (\ref{i6}) resp. (\ref{w5+}) to be well
defined, the functions $r$, $\vc{U}$ must belong at least to the
class \bFormula{b1+} r -\overline\vr\in C_{\rm weak}([0,T]; L^2+
L^\gamma (\Omega)),\;  \eF
\bFormula{b2} \vc{U} \in \ L^2(0,T; W^{1,2}(\Omega;R^3)).
%
%
\eF
\item
A short inspection (\ref{i7}) resp. (\ref{w6+}) implies that the integrals are well-defined
if, at least,
\bFormula{b3+}
\partial_t \vc{U} \in L^2(0,T; L^{3}\cap L^{6 \gamma/  (5 \gamma - 6)}(\Omega,
R^3))+ L^1(0,T; L^{4/3}\cap L^{2 \gamma/  (\gamma - 1)}(\Omega,
R^3)),
\eF
\bFormula{b5+} \Grad \vc{U} \in  L^\infty(0,T; L^{6}\cap
L^{3 \gamma/ (2 \gamma - 3)}(\Omega, R^{3 \times 3})) + L^2(0,T; L^{12/7}
\cap L^{6 \gamma/ (4 \gamma - 3)}(\Omega, R^{3 \times 3}))
\eF
\[
 + L^1(0,T;
L^\infty(\Omega; R^3)), \]
\bFormula{b4} {\rm div} \vc{U} \in L^1(0,T; L^\infty(\Omega)),
\eF
\item
The function $r$ must be bounded below away from zero, and
\bFormula{b6}
\partial_t H'(r) \in L^1(0,T; L^{\gamma/ (\gamma -
1)}\cap  L^2(\Omega)),
\eF
\bFormula{b7} \Grad H'(r) \in L^2(0,T; L^{3}\cap L^{6 \gamma/  (5
\gamma - 6)}(\Omega, R^3))+ L^1(0,T; L^{4/3}\cap L^{2 \gamma/
(\gamma - 1)}(\Omega, R^3)). \eF

\item Finally, the vector field $\vc U$ has to satisfy
\begin{equation}\label{b8}
\begin{array}{c}
\vc U|_{\partial\Omega}=0 \;\mbox{in the case of boundary
conditions (\ref{i4})},
\\ \\
\vc U\cdot\vc n|_{\partial\Omega}=0 \;\mbox{in the case of
boundary conditions (\ref{p6})}.
\end{array}
\end{equation}
\end{itemize}

Consequently, Theorem \ref{Tw1} is valid even if we replace the
hypotheses on smoothness and integrability of the test functions $(r,\vc U)$ by
weaker hypotheses, namely  (\ref{b1+}--\ref{b8}).

In particular, $r$, $\vc{U}$ may be another (strong) solution
emanating from the same initial data $\vr_0$, $\vu_0$. Specific
examples will be discussed in the forthcoming section.

\section{Applications}
\label{a}

In this section, we show how Theorem \ref{Tw1} can be applied in order to establish weak-strong uniqueness property for the compressible Navier-Stokes
system in the class of finite energy weak solutions in bounded and unbounded domains. Other applications can be found in \cite{FENOSU}.

\subsection{Weak-strong uniqueness on bounded domains}
\subsubsection{No-slip boundary conditions}
To begin, observe that \emph{any} finite energy weak solution $\vr$, $\vu$ of the compressible Navier-Stokes system (\ref{i1} - \ref{i4}) in $(0,T)
\times \Omega$, where $\Omega$ is a bounded domain, belongs to the class
\[
\vr \in C_{\rm weak}([0,T]; L^\gamma(\Omega)),\
\vr \vu \in C_{\rm weak}([0,T]; L^{2\gamma/(\gamma + 1)}(\Omega;R^3)), \
\vu \in L^2(0,T;W^{1,2}_0(\Omega;R^3)),
\]
and, by virtue of the energy inequality (\ref{w4}),
\[
 p(\vr) \in L^\infty(0,T; L^1(\Omega)).
\]

Moreover, it is easy to check that
\bFormula{a1}
H(\vr) - H'(r)(\vr - r) - H(r) \geq c(r)
\left\{ \begin{array}{l} (\vr - r)^2 \ \mbox{for}\ r/2 < \vr < 2 r,
\\ \\ (1 + \vr^\gamma) \ \mbox{otherwise}
\end{array} \right.,
\eF
where $c(r)$ is uniformly bounded for $r$ belonging to compact sets in $(0, \infty)$.

Finally, note that, since the total mass is a conserved quantity on a bounded domain, we can take $\Ov{\vr}$ in (\ref{i5a}) so that
\[
\intO{ (\vr - \Ov{\vr}) } = 0.
\]

The rather obvious leading idea of the proof of weak-strong uniqueness is to take $r = \tilde \vr$, $\vc{U} = \tilde \vu$ in the relative energy
inequality (\ref{i6}), where $\tilde \vr$, $\tilde \vu$ is a (hypothetical) regular solution, originating from the same initial data. The following
formal computations will require certain smoothness of $\tilde \vr$, $\tilde \vu$ specified in the concluding theorem. Moreover, we assume that
$\tilde \vr$ is bounded below away from zero on the whole compact time interval $[0,T]$.

Our goal is to examine all terms in the remainder (\ref{i7}) and to show they can be ``absorbed'' by the left-hand side of (\ref{i6}) by means of a
Gronwall type argument.

\begin{enumerate}
\item We rewrite
\[
\intO{  \vr \Big( \partial_t \tilde \vc{u} + \vu \cdot \Grad \tilde \vc{u} \Big) \cdot
(\tilde \vc{u} - \vu )} = \intO{  \vr \Big( \partial_t \tilde \vc{u} + \tilde \vc{u} \cdot \Grad \tilde \vc{u} \Big) \cdot
(\tilde \vc{u} - \vu )} + \intO{  \vr (\vc{u} - \tilde \vc{u}) \cdot \Grad \tilde \vc{u}  \cdot
(\tilde \vc{u} - \vu )}.
\]

Seeing that
\[
\partial_t \tilde \vc{u} + \tilde \vc{u} \cdot \Grad \tilde \vc{u} =
\frac{1}{\tilde \vr} \Div \tn{S} (\Grad \tilde \vu) + \vc{f} - \Grad H'(\tilde \vr) ,
\]
we go back to (\ref{i7}) to obtain
\[
{\cal R}(\vr, \vu, \tilde \vr, \tilde \vu) = \intO{  \vr (\vc{u} - \tilde \vc{u}) \cdot \Grad \tilde \vc{u}  \cdot
(\tilde \vc{u} - \vu )} +  \intO{ \frac{1}{\tilde \vr}\left( \vr  - \tilde \vr \right) \Div \tn{S} (\Grad \tilde \vu) \cdot (\tilde \vu - \vu) }
\]
\[
+ \intO{(\tilde \vr - \vr) \Big( \partial_t H'(\tilde \vr) + \Grad H'(\tilde \vr)
\cdot \tilde \vu \Big)} - \intO{ \Div \tilde \vu \Big( p(\vr) - p(\tilde \vr) \Big) }.
\]

\item
Computing
\[
(\tilde \vr - \vr) \Big( \partial_t H'(\tilde \vr) + \Grad H'(\tilde \vr)
\cdot \tilde \vu \Big) = - \Div \tilde \vu (\vr - \tilde \vr) p'(\tilde \vr),
\]
we may infer that
\[
\intO{(\tilde \vr - \vr) \Big( \partial_t H'(\tilde \vr) + \Grad H'(\tilde \vr)
\cdot \tilde \vu \Big)} - \intO{ \Div \tilde \vu \Big( p(\vr) - p(\tilde \vr) \Big) }
\]
\[
= - \intO{ \Div \tilde \vu \Big( p(\vr) - p'(\tilde \vr) (\vr - \tilde \vr) - p(\tilde \vr) \Big) };
\]
whence
\bFormula{a2}
{\cal R} (\vr, \vu, \tilde \vr, \tilde \vu) = \intO{  \vr (\vc{u} - \tilde \vc{u}) \cdot \Grad \tilde \vc{u}  \cdot
(\tilde \vc{u} - \vu )} - \intO{ \Div \tilde \vu \Big( p(\vr) - p'(\tilde \vr) (\vr - \tilde \vr) - p(\tilde \vr) \Big) }
\eF
\[
+\intO{ \frac{1}{\tilde \vr}\left( \vr  - \tilde \vr \right) \Div \tn{S} (\Grad \tilde \vu) \cdot (\tilde \vu - \vu) }.
\]

\item
In view of (\ref{a1}), we have
\bFormula{a3}
\left|
 \intO{  \vr (\vc{u} - \tilde \vc{u}) \cdot \Grad \tilde \vc{u}  \cdot
(\tilde \vc{u} - \vu )} - \intO{ \Div \tilde \vu \Big( p(\vr) - p'(\tilde \vr) (\vr - \tilde \vr) - p(\tilde \vr) \Big) } \right|
\eF
\[
\leq c \| \Grad \tilde \vu \|_{L^\infty(\Omega;R^3)} {\cal E}
\Big( [\vr, \vu] \Big| [\tilde \vr, \tilde \vu ] \Big),
\]
provided
\bFormula{a4}
0 < \inf_{[0,T] \times \Ov{\Omega}} \tilde \vr \leq \tilde \vr (t,x)\leq
\sup_{[0,T] \times \Ov{\Omega}} \tilde \vr < \infty.
\eF

\item Finally, we write
\[
\intO{ \frac{1}{\tilde \vr}\left( \vr  - \tilde \vr \right) \Div \tn{S} (\Grad \tilde \vu) \cdot (\tilde \vu - \vu) }
\]
\[
= \int_{ \{ \tilde \vr / 2 < \vr < 2 \tilde \vr \}}  \frac{1}{\tilde \vr}\left( \vr  - \tilde \vr \right) \Div \tn{S} (\Grad \tilde \vu) \cdot
(\tilde \vu - \vu) \ \dx
\]
\[
+ \int_{ \{ 0 \leq \vr \leq \tilde \vr/2 \}}  \frac{1}{\tilde \vr}\left( \vr  - \tilde \vr \right) \Div \tn{S} (\Grad \tilde \vu) \cdot (\tilde
\vu - \vu) \ \dx + \int_{ \{  \vr \geq 2 \tilde \vr  \}}  \frac{1}{\tilde \vr}\left( \vr  - \tilde \vr \right) \Div \tn{S} (\Grad \tilde \vu)
\cdot (\tilde \vu - \vu) \ \dx,
\]
where, by virtue of H\" older's inequality,
\bFormula{a5}
\left|
\int_{ \{ \tilde \vr / 2 < \vr < 2 \tilde \vr \}}  \frac{1}{\tilde \vr}\left( \vr  - \tilde \vr \right) \Div \tn{S} (\Grad \tilde \vu) \cdot
(\tilde \vu - \vu) \ \dx
\right|
\eF
\[
\leq c(\delta) \left\| \frac{1}{\tilde \vr}  \Div \tn{S}
(\Grad \tilde \vu) \right\|_{L^3(\Omega;R^3)}^2  \int_{ \{ \tilde \vr / 2 < \vr < 2 \tilde \vr \}}
(\vr - \tilde \vr)^2 \ \dx + \delta \| \tilde \vu - \vu \|^2_{L^6(\Omega;R^3)}
\]
for any $\delta > 0$.

Furthermore, in accordance with (\ref{a1}), we get \bFormula{a6}
\int_{ \{ \tilde \vr / 2 < \vr < 2 \tilde \vr \}} (\vr - \tilde
\vr)^2 \ \dx \leq c  {\cal E} \Big( [\vr, \vu] \Big| [\tilde \vr,
\tilde \vu ] \Big), \eF while, by virtue of Sobolev's inequality
and Korn-type inequality (see e.g. Dain \cite{DAIN})
\begin{equation}\label{korn}
\|\vc z\|_{1,2}\le c\|\tn S(\Grad \vc z)\|_{L^2(\Omega;R^{3\times
3})},\; \vc z\in W^{1,2}(\Omega;R^3),
\end{equation}
we have
 \bFormula{a7} \| \tilde \vu - \vu \|^2_{L^6(\Omega;R^3)} \leq
c \| \Grad \vu - \Grad \tilde \vu \|^2_{L^2(\Omega; R^{3 \times
3})}\le c\|\tn S(\Grad \vc u-\Grad\tilde\vc
u)\|^2_{L^2(\Omega;R^{3\times 3})} . \eF
Therefore,
\[
\left| \int_{ \{0\le \vr \leq \tilde \vr/2\}} \frac{1}{\tilde
\vr}\left( \vr  - \tilde \vr \right) \Div \tn{S} (\Grad \tilde
\vu) \cdot (\tilde \vu - \vu) \ \dx \right|\le \]
\[
\leq c(\delta) \left\| \frac{1}{\tilde \vr}\Div \tn{S} (\Grad
\tilde \vu) \right\|_{L^3(\Omega;R^3)}^2 {\cal E} \Big( [\vr, \vu]
\Big| [\tilde \vr, \tilde \vu ] \Big) + \delta  \|\tn S (\Grad \vu
- \Grad \tilde \vu )\|^2_{L^2(\Omega; R^{3 \times 3})}
\]
for any $\delta > 0$.

 Next we realize that
$$
{\cal E}(\vr,\vt|\tilde\vr,\tilde\vt)\in L^\infty(0,T)
$$
and that
$$
\|\vr\|_{L^\gamma(\{\vr>2\overline\vr\})}\le c \Big[{\cal
E}(\vr,\vt| \tilde \vr,\tilde\vt)\Big]^{1/\gamma}, \quad
\|\vr^{\gamma/2}\|_{L^2(\{\vr>2\overline\gamma\})}\le c\Big[ {\cal
E}(\vr,\vt| \tilde \vr,\tilde\vt)\Big]^{1/2}.
$$
Using these facts, we  deduce
\bFormula{a9} \left| \int_{ \{ \vr \geq 2\tilde \vr\}}
\frac{1}{\tilde \vr}\left( \vr  - \tilde \vr \right) \Div \tn{S}
(\Grad \tilde \vu) \cdot (\tilde \vu - \vu) \ \dx \right|\le \eF
\[
 \int_{\{ \vr \ge 2\tilde\vr \} }\left(
\left| \frac{\vr - \tilde \vr}{\vr\tilde\vr}\right|{\rm
max}\{\vr,\vr^{\gamma/2}\}\left| \Div \tn{S}(\Grad \tilde
\vu)\right| \,\left|(\tilde \vu - \vu ) \right|\right)(\tau,\cdot)
\ \dx\le
\]
\[
c\|\tn S(\Grad \vu - \Grad \tilde \vu )\|_{L^2(\Omega;
R^{3\times3})} \| \Div \tn{S}(\Grad \tilde \vu) \|_{L^q\cap
L^3(\Omega; R^3)} \Big[{\cal E}(\vr,\vt| \tilde
\vr,\tilde\vt)\Big]^{1/2}\le
\]
\[
\le \delta \| \tn S(\Grad\vu -\Grad \tilde \vu \|_{L^2(\Omega;
R^3))}^2 + c(\delta) \| \Div \tn{S}(\Grad \tilde \vu)
\|^2_{L^q\cap L^3(\Omega; R^3)}\; {\cal  E}(\vr,\vt| \tilde
\vr,\tilde\vt),\; q = \frac{6 \gamma}{5 \gamma - 6}.
\]


\end{enumerate}

Summing up relations (\ref{a2} - \ref{a9}) we conclude that the relative entropy inequality, applied to $r = \tilde \vr$, $\vc{U} = \tilde \vu$,
yields the desired conclusion
\bFormula{a10}
{\cal E} \Big( [\vr, \vu] \Big| [\tilde \vr, \tilde \vu] \Big) (\tau) \leq
\int_0^\tau h(t) {\cal E} \Big( [\vr, \vu] \Big| [\tilde \vr, \tilde \vu] \Big) (t) \ \dt, \ \mbox{with}\ h \in L^1(0,T),
\eF
provided $\tilde \vr$ satisfies (\ref{a4}), and
\bFormula{a11}
\Grad \tilde \vu \in L^1 (0,T; L^\infty (\Omega; R^{3 \times 3})) \cap
L^2(0,T; L^2(\Omega; R^{3 \times 3})),\
\Div \tn{S} (\Grad \tilde \vu) \in L^2(0,T; L^3 \cap L^q(\Omega;R^3)),
\eF
with
\[
q = \frac{6 \gamma}{5 \gamma - 6}.
\]

We have shown the following result:

\bTheorem{a1}
Let $\Omega \subset R^3$ be a bounded Lipschitz domain, let the pressure $p$ satisfy hypothesis (\ref{w1}), and let
\[
\vc{f} \in L^1(0,T; L^{2\gamma/(\gamma - 1)}(\Omega;R^3)).
\]
Assume that $\vr$, $\vu$ is a finite energy weak solution to the Navier-Stokes system (\ref{i1} - \ref{i4}) in $(0,T) \times \Omega$, specified in
Section \ref{fws}. Let $\tilde \vr$, $\tilde \vu$ be a (strong) solution of the same problem belonging to the class
\[
0 < \inf_{(0,T) \times \Omega} \tilde \vr \leq \tilde \vr (t,x) \leq
\sup_{(0,T) \times \Omega} \tilde \vr < \infty,
\]
\[
\Grad \tilde \vr \in L^2(0,T; L^q(\Omega; R^3)),\
\Grad^2 \tilde \vu \in L^2(0,T; L^q(\Omega; R^{3 \times 3 \times 3})),
\ q > \max \left\{ 3 ; \frac{3}{\gamma - 1} \right\},
\]
emanating from the same initial data.

Then
\[
\vr = \tilde \vr , \ \vu = \tilde \vu \ \mbox{in}\ (0,T) \times \Omega.
\]
\eT

\medskip

{\bf Remark \ref{a}.1} {\it We need $\Omega$ to be at least
Lipschitz to guarantee the $W^{1,p}$ extension property, with the
associated embedding relations}.

\medskip

{\bf Remark \ref{a}.2} {\it
The reader will have noticed that the regularity properties required for $\tilde \vr$, $\tilde \vu$ in Theorem \ref{Ta1} are in fact \emph{stronger}
than (\ref{a11}). The reason is that all integrands appearing in the relative energy inequality (\ref{i6}) must be well defined.}

\medskip

{\bf Remark \ref{a}.3} {\it
\emph{Existence of finite energy weak solutions} was shown in \cite{FNP1} for
general (finite energy) data and without any restriction on imposed on smoothness of $\partial \Omega$.}

\medskip

{\bf Remark \ref{a}.4} {\it \emph{Local-in-time existence of
strong solutions} belonging to the regularity class specified in
Theorem \ref{Ta1} was proved by Sun, Wang, and Zhang
\cite{SuWaZh}, under natural restrictions imposed on the initial
data.}

\subsubsection{Navier boundary conditions with friction}

Theorem \ref{Ta1} holds in the case of Navier's boundary condition
(\ref{p6}).
The proof remains basically without changes; the standard Korn-type
inequality (\ref{korn})  has to be replaced by a more
sophisticated one, namely
\begin{equation}\label{korng}
\begin{array}{c}
\|\vc v\|^2_{W^{1,2}(\Omega,R^3)}\le c(M,K,p)\Big(\|\tn S(\Grad\vc
v)\|^2_{L^2(\Omega,R^{3\times 3})}+\| R \vc
v^2\|_{L^1(\Omega)}\Big)\\ \\
\mbox{for any $\vc v\in W^{1,2}(\Omega;R^3)$, $R\ge 0$,
$M\le\int_\Omega R {\rm d}x$, $\|R\|_{L^p(\Omega)}\le K$},
\end{array}
\end{equation}
where $M,K>0$, $p>1$ (see \cite[Theorem 10.17]{FEINOV}). It is
employed in estimate (\ref{a7}) with $\vc v=\vc u-\tilde\vc u$ and
$R=\vr$.

\subsection{Weak strong uniqueness on unbounded domains}
\subsubsection{No-slip boundary conditions}

If the Navier-Stokes system is considered on an unbounded domain $\Omega$, the far-field behavior (\ref{i4a}) must be specified. Here, we assume that
$\Ov{\vr} > 0$ so that the density $\tilde \vr$ of the (hypothetical) strong solution may be bounded below away from zero. Moreover, the finite energy
weak solutions necessarily belong to the class:
\bFormula{ar1}
\vr - \Ov{\vr} \in L^\infty(0,T;L^2 + L^\gamma(\Omega)),\
p(\vr) - p(\Ov{\vr}) \in L^\infty(0,T; L^2 + L^1 (\Omega)),
\eF
\bFormula{ar2}
\vu \in L^2(0,T;W^{1,2}_0(\Omega;R^3)),\ \vr \vu \in
L^\infty(0,T;L^2 + L^{2\gamma/(\gamma + 1)}(\Omega;R^3)).
\eF

An appropriate modification of Theorem \ref{Ta1} for unbounded
domains reads: \bTheorem{a2} Let $\Omega \subset R^3$ be an
unbounded domain with a uniformly Lipschitz boundary, let the
pressure $p$ satisfy hypothesis (\ref{w1}), and let
\[
\vc{f} \in L^1(0,T; L^{1} \cap L^\infty (\Omega;R^3)).
\]
Assume that $\vr$, $\vu$ is a finite energy weak solution to the Navier-Stokes system (\ref{i1} - \ref{i4}) in $(0,T) \times \Omega$, specified in
Section \ref{fws}, satisfying the far-field boundary
conditions (\ref{i4a}), with $\Ov{\vr} > 0$. Let $\tilde \vr$, $\tilde \vu$ be a (strong) solution of the same problem belonging to the class
\[
0 < \inf_{(0,T) \times \Omega} \tilde \vr \leq \tilde \vr (t,x) \leq
\sup_{(0,T) \times \Omega} \tilde \vr < \infty,
\]
\[
\Grad \tilde \vr \in L^2(0,T; L^2 \cap L^q(\Omega; R^3)),\
\Grad^2 \tilde \vu \in L^2(0,T; L^2 \cap L^q(\Omega; R^{3 \times 3 \times 3})),
\ q > \max \left\{ 3 ; \frac{3}{\gamma - 1} \right\},
\]
emanating from the same initial data, and satisfying the energy inequality (\ref{i8}).

Then
\[
\vr = \tilde \vr , \ \vu = \tilde \vu \ \mbox{in}\ (0,T) \times \Omega.
\]
\eT

\medskip

{\bf Remark \ref{a}.5} {\it The uniformly Lipschitz boundary
$\partial \Omega$ guarantees the $W^{1,p}$-
extension property as well as validity of Korn's inequality
(\ref{korn}). }

\medskip

{\bf Remark \ref{a}.6} {\it Since the strong solution satisfies the energy
(in)equality (\ref{i8}), it automatically belongs to the regularity class
(\ref{ar1}), (\ref{ar2}). }

\medskip

{\bf Remark \ref{a}.7} {\it \emph{Existence of finite energy weak solutions}
for certain classes of unbounded domains was shown in  \cite{NOST4}, see also Lions \cite{LI4}.}

\medskip

{\bf Remark \ref{a}.8} {\it The reader may consult the nowadays classical papers by Matsumura and Nishida
\cite{MANI1}, \cite{MANI} for the existence of strong solutions,
more recent results can be found in Cho, Choe and Kim \cite{ChoChoeKim}, and in the references cited therein.}

\subsubsection{Navier boundary conditions}

Theorem \ref{Ta2} remains valid also for the Navier boundary conditions. We have however suppose that
on the considered unbounded domain a sort of Korn type inequality
holds, for example
\bFormula{korngu} \| \vc{v} \|_{W^{1,2}(\Omega;R^3)}^2 \leq c(|V|)
\left( \| \tn S (\Grad \vc{v} ) \|^2_{L^2(\Omega;R^3)} +
\int_{\Omega \setminus V} |\vc{v}|^2 \ \dx  \right), \eF
\[
\mbox{for any}\ \vc{v} \in W^{1,2}(\Omega),  \ |V| < \infty.
\]
Such inequality is known to hold  in a half space, an exterior
domain, a cylinder, a plane slab, to name only a few.

Since
$$
\Big|\{|\vr-\overline\vr|\ge \overline\vr/2\}|<\infty,
$$
inequality (\ref{korngu}) implies the validity of (\ref{korng})
with  $\vc v=\vc u-\tilde\vc u$ and $R=\vr$. This inequality has
to replace the standard Korn's inequality (\ref{korn}) in estimate
(\ref{a7}). Other arguments in the proof remain without changes.

\def\ocirc#1{\ifmmode\setbox0=\hbox{$#1$}\dimen0=\ht0 \advance\dimen0
  by1pt\rlap{\hbox to\wd0{\hss\raise\dimen0
  \hbox{\hskip.2em$\scriptscriptstyle\circ$}\hss}}#1\else {\accent"17 #1}\fi}



\end{document}